\documentclass[12pt]{article}
\usepackage{amsmath,latexsym,amssymb,graphicx,epsf,amsfonts,amsthm}
\usepackage{epsf}
\usepackage{color}
\numberwithin{equation}{section}

\newtheorem*{ass}{Assumption}

\begin{document}
\hoffset = -1truecm \voffset = -2truecm
\title{A Sparse Semi-Blind Source Identification Method and Its Application to
Raman Spectroscopy for Explosives Detection}
\author{Yuanchang Sun \thanks{Department of Mathematics,
University of California at Irvine, Irvine, CA 92697, USA.} and Jack
Xin$^{*}$}
\date{}
\maketitle
\begin{abstract}
Rapid and reliable detection and
identification of unknown chemical substances is critical to homeland security.
It is challenging to identify chemical components
from a wide range of explosives.  There are two key steps involved.
One is a nondestructive and informative
spectroscopic technique for data acquisition. The other is an associated library of
reference features along with a computational method for feature matching and
meaningful detection within or beyond the library.

Recently several experimental techniques based on Raman scattering have been developed to
perform standoff detection and identification of explosives, and they prove to be
successful under certain idealized conditions.  However data
analysis is limited to standard least squares method assuming the complete knowledge of the chemical components.
In this paper, we develop a new iterative method to identify unknown substances from
mixture samples of Raman spectroscopy. In the first step, a constrained least squares method
decomposes the data into a sum of linear combination of the known components
and a non-negative residual. In the second step, a sparse and convex blind source separation method
extracts components geometrically from the residuals. Verification based on the library templates or expert
knowledge helps to confirm these components. If necessary, the confirmed meaningful components are fed back into
step one to refine the residual and then step two extracts possibly more hidden components. The
two steps may be iterated until no more components can be identified.
 We illustrate the proposed method in
processing a set of the so called swept wavelength optical resonant Raman spectroscopy experimental data
by a satisfactory blind extraction
of a priori unknown chemical explosives from mixture samples.
\end{abstract}
\thispagestyle{empty}
\newpage
\setcounter{equation}{0} \setcounter{page}{1}
\section{Introduction}
A critical problem in homeland security is reliable and rapid identification
of unknown chemical and biological substances in the explosives.
Due to the harmful environment caused by the release of the explosive chemicals,
non-destructive spectroscopic techniques are typically used
to record the optical spectrum without interfering with the samples.
Ideally, a standoff detection is performed to acquire the spectral information.
Then a search and matching procedure through a prepared spectral database
will be carried out for identification. However such an approach would be unsuccessful if
the explosives contained chemical compounds outside of the database, which is highly likely
as hidden explosives are often unknown a-priori. In general, the samples may
involve multiple unknown substances besides impurities.
Conventional analysis routines are mostly based on least squares fitting whose
residuals could remain mysterious. Further analysis calls for the development of
blind identification methods to extract major components from the residuals.

\medskip
Various recent experimental techniques
(see \cite{Comanescu_08,Gaft_07,Grun_07} and the references therein) can
identify pure chemicals with notable
success. These methods are mainly based on Raman spectroscopy,
a spectroscopic technique to study the
chemical composition of the samples \cite{long_02, Turrel_72}.
Combined with other spectroscopic techniques such as Ultraviolet and Infrared spectroscopy,
Raman spectroscopy has been widely used in materials science, biosciences,
geosciences (gemology), forensic sciences, nano-technology,
and pharmaceutical chemistry \cite{Comanescu_08, Faulds_04, Froud_03, Grun_07, Miller_01, Not_03, Sas_06, Yao_07}.
For example, the Swept Wavelength Optical Resonant Raman Detector (SWOrRD) at
the Naval Research Laboratory developed in 2009 can generate two dimensional
spectral maps of biological agents and chemical substances.
The resultant two dimensional signatures contain much more information than the
single illumination wavelengths, which may result in a
greater likelihood of successful identification even in
complex mixtures \cite {Comanescu_08,Grun_07}.
Raman spectroscopy is based on Raman scattering, an inelastic scattering process
that shifts the frequencis of the incident photons. During the interaction of
the incident light with a molecule, a scattered light of lower (Stokes Raman) or
higher (anti-Stokes Raman) energy is emitted, allowing the measurements of the molecule's
 vibrational modes. The appeals of this technique are non-destructiveness and
fast sensing capability. It has become a promising tool for standoff distance detection
for explosives at airports among other transportation centers.

\medskip
In general, Raman spectra of a sample are composed of many substances, and
they must be identified by an analysis software.
If one has the complete knowledge of the kinds of substances in the sample,
the least squares fitting can be used to retrieve their concentrations (or volumes)
by a linear combination of the known spectra on a template.
In most practical situations, one may have to identify the substances and quantify
their concentrations at the same time.
This becomes a blind source separation (BSS) problem, or
recovering pure signal sources from their mixtures
without a detailed knowledge of the mixing process.  There have been several studies on
the BSS of Raman spectra  \cite{Miron_11, Rowland_11} based
on independent component analysis (ICA, \cite{Hyv}) and nonnegative matrix
factorization (NMF, \cite{Lee}).  However these methods
are non-convex and too general to be robust and reliable
in real-world applications. The independence hypothesis of ICA
does not hold if chemicals share some common structures and
have correlated Raman spectral line shapes.  Moreover, these existing approaches do not
address how to identify unknown substances from the fitting residuals
when partial knowledge of the source signals is available.
Such a semi-blind problem is more often encountered in applications and
of great importance to practitioners.

\medskip
In this paper, we shall develop a convex semi-blind source separation method
based on sparsity of the source spectra.
We are concerned with the regime
where the sample contains some known and some unknown components.
In other words, we have knowledge of some of the components and
their concentrations, which is the case of the SWOrRD data..
We further assume that the upper bound of the concentrations
of the known substances is available, say
from experiments or prior knowledge, as is the case of SWOrRD data.
Though our method here is designed for Raman spectra,
it is applicable to spectra with similar line shapes,
such as nuclear magnetic resonance (NMR) data, see \cite{S_nmr_1,S_nmr_2} where
more general sparseness source conditions and post-processing methods have been studied.
A similar semi-blind source
identification approach has been developed for
mysterious species of the atmospheric gas mixtures \cite{Sun_Xin_Doas}.
\bigskip

Our method in the paper may provide an assistive tool for practitioners to
analyze the hidden structures from residuals of fitting complicated mixtures.
The method contains two major steps.
The first step decomposes Raman data $X$ into
a linear combination of the reference spectra (known components) plus a remainder,
or  $X = A\,S + R$, where the columns of matrix $A$ are known reference spectra,
the matrix $S$ contains the non-negative concentrations and has known upper bound.
The first step is carried out by solving a constrained
least squares problem.  The constraints on $S$ help to maintain the
nonnegativity of remainder matrix $R$.
The second step performs a nonnegative blind source separation
of the remainder matrix $R$ to extract the unknown components.
We show that proper sparsity of source signals
reduces the general non-convex problem to constrained convex programming
permitting solutions with better mathematical properties.
Sparse solutions to convex objectives are achieved
through $\ell_1$ norm minimization and a fast iterative method (the linearized Bregman).
The two steps may be iterated. If some of the components from step two are confirmed as
chemically meaningful, they are fed back to step one to refine the residual for
further extraction of hidden components in step 2.

\medskip
The paper is organized as follows.  In section 2,
we introduce the source sparseness condition and our method.
In section 3, we illustrate the proposed method in processing a
set of SWorRD experimental data, and
show satisfactory numerical results.
Conclusion and future works are discussed in section 4.

\medskip
This work was partially supported by NSF-ADT grant DMS-0911277.
The authors thank Naval Research Laboratory for the experimental Raman data and Dr. J. Grun for helpful discussions.

\section{The method}
In this section, we shall present our method for chemicals detection in Raman spectroscopy.
Based on the assumption that part of the substances in the chemical mixture and
the upper bounds of their concentrations
are known, the first step of the method fits these known chemicals to the mixtures.
This step solves a constrained least squares optimization.
\subsection{Constrained least squares}
The following linear model is used for Raman spectra of mixtures
\begin{equation}
\label{Ls}
X =  A\,S + R,
\end{equation}
where the columns of matrix $X$ are the measured Raman data,
the columns of $A$ contain reference spectra of the known chemicals,
and those of $S$ matrix contains their concentrations.
Matrix $R$ is the fitting residual which might contain hidden spectral structures,
the instrument noise, etc.
Matrix $S$ is nonnegative on physical ground that
its entries represent concentrations.
Furthermore, the upper bound of $S$ is known {\it a priori}.
Then for the estimation of $S$, we minimize the following constrained objective function:
\begin{equation}
\label{convexModel}
\min_{S} \|X-A\,S\|^2_2,\;\; \mathrm{s.t.}\;\;0\leq S \leq c,
\end{equation}
where the vector $c$ contains the upper bounds of the concentrations of the
known substances in the sample.
For the Raman data in our numerical experiments,
we found that the linear constraint in (\ref{convexModel}) helps to
maintain a nonnegative remainder $R = X-A\,S$.  In the residual $R$, there might
be spectral structures (one or many) of chemicals buried in noise, or just random noise.
 In either case, we factorize the residuals in a blind fashion
due to the lack of knowledge of the hidden components.
Conventional blind source separation methods such as NMF and ICA
are non-convex optimization methods.  These methods are for general purpose, yet
often unreliable in real-world applications due to non-convexity or sensitivity of
their working assumptions. For our problem, we show that proper sparsity of
source signals reduces the general non-convex problem to constrained convex programming
permitting solutions with better mathematical properties.
Sparse solutions to convex objectives are achieved through $\ell_1$ minimization.
\subsection{Convex blind source separation}
The remainder matrix $R$ is factorized as:
\begin{equation}
\label{modelR}
R = W\,M,
\end{equation}
where the columns of the matrix $W$ are the identified substances,
$M$ is their concentrations in the sample.  All the matrices are nonnegative.
For the purpose of illustration, we shall call $W \in \mathbb{R}^{p\times n}$ the source matrix, $M\in \mathbb{R}^{n\times m}$ the mixing matrix.  The dimensions of the matrices are expressed in terms of three numbers:
(1) $p$ the number of available samples, (2) $m$ the number of mixture
signals, and (3) $n$ the number of source signals.  For the Raman data we considered in the paper, there are more mixtures than sources, i.e., $m\geq n$.  The goal is to recover $W$ and $M$ for a given $R$.  This is also known as an NMF problem. Various methods have been developed to solve BSS problems by exploiting the natures of source signals.  For example, independent component analysis (ICA) recovers statistically independent signals.  However, the independence should not be assumed on the Raman spectra of the chemical substances when they share common structural features (the line shapes of their Raman spectra are similar).  A better working assumption for the data is a so-called partial sparseness condition.  Namely, the source signals are only required to be non-overlapping at some locations of acquisition variable.  This sparseness condition was first known in the 1990s \cite{Boardman_93,Winter_99} in the study of blind hyper-spectral unmixing of remote sensing, where the source condition is called pixel purity assumption (PPA) \cite{Chang_07}.   In 2005,  Naanaa and Nuzillard \cite{NN05} used this assumption to separate the signals in nuclear magnetic resonance spectroscopy.
In fact, this condition is well suited to many chemical substances including those studied in this paper.  Such a sparseness condition leads to a dramatic mathematical simplification of a general non-negative matrix factorization problem (\ref{modelR}) which is non-convex.  Geometrically speaking, the problem of finding the mixing matrix $M$ reduces to
the identification of a minimal cone containing the rows of matrix $R$.  The latter can be achieved by linear programming.  In the context of hyper-spectral unmixing, the resulting geometric (cone) method is the so called N-findr \cite{Winter_99}, and is now a benchmark in hyper-spectral unmixing.  Next we shall review the essentials of the partial spareness condition and the geometric cone method.

\medskip
Simply speaking, the key sparseness assumption on the
source signals is that each source has a stand alone peak at some
location of acquisition variable where the other sources are
zero.  More precisely, the source matrix $W\geq 0$ is
assumed to satisfy the following condition:
\begin{ass}
For each $j\in\{1,2,\dots,n \}$ there exists an $i_{j}\in
\{1,2,\dots,p\}$ such that $w_{i_j,j}>0$ and $w_{i_j,k}=0,\;\;\;\;
(k=1,\dots,i-1,i+1,\dots,n).$
\end{ass}
Eq. (\ref{modelR}) can be rewritten in terms of rows as
\begin{equation}
\label{LinComb}
R^{i} = \sum^{n}_{k=1}w_{i,k}M^{k}\;\; i = 1,\dots,p\;,
\end{equation} where $ R^i$ denote the $i$th row of $R$, and $M^k$ the $k$th row of $M$.
The source assumption implies that $\displaystyle R^{i_j} = w_{i_j,j}M^j\;\; j = 1,\dots,n\;\; $ or $M^{j} = \frac{1}{w_{i_j,j}}R^{i_j}$.  Hence
Eq. (\ref{LinComb}) is rewritten as
\begin{equation}
\label{NNlinComb} R^{i} = \sum^{n}_{k = 1} \frac{w_{i,k}}{w_{i_k,k}} M^{i_k}\;,
\end{equation}
which says that every row of $R$ is a nonnegative linear combination of the rows of $M$.  The identification of $ M $'s rows is equivalent to identifying a convex cone of a finite collection of vectors.  The cone encloses the data rows in matrix $R$, and is the smallest of such cones.  Such a minimal enclosing convex cone can be found by linear programming methods.  Mathematically, the following optimization problems are suggested to estimate the mixing matrix
{\allowdisplaybreaks
\begin{equation}
\label{LPNP} \min\; \mathrm{score} = \|\sum^{p}_{i = 1, i\neq l}R^{i} \lambda_i - R^{l}
 \|_2\;, \mathrm{s.t.} \;\; \lambda_i \geq 0\;,  l = 1,\dots,p\;.
\end{equation}}
A score is associated with each row of $R$.  A row with a low score is unlikely to be a row of $M$ because this
row is roughly a nonnegative linear combination of the other rows of $R$.  On the other hand, a high score means that the
corresponding row is far from being a nonnegative linear combination of other rows.
The $n$ rows from $R$ with highest scores are selected to form $M$, the mixing matrix.
\subsection{Sparse source recovery}
With the recovery of $M$, we solve for $W$ next.
An existing method \cite{NN05} is to directly compute $W = R\,M^{+}$, where $M^{+}$ is the
pseudo-inverse of $M$.  It is however sensitive to noise,
and tends to introduce errors and artifacts of negative values.
Although a nonnegative least squares can produce a nonnegative $W$,
spurious peaks might be introduced in the results.
To benefit from the sparseness source condition in the first step,
a more reliable method is to solve a nonnegative $\ell_1$ optimization.
Although the source signals (columns of $W$) are not sparse,
the rows of $W$ possess sparsity.  Hence, we seek the sparsest solution for
each row $W^i$ of $W$.
\begin{equation}
\label{Lzero}
 \min \|W^i \|_0\quad \mathrm{subject}\;\mathrm{to}\;  W^i\,M = R^i,\;\; W^i \geq 0.
\end{equation}
Here $\|\cdot\|_0$ ( 0-norm ) represents the number of nonzeros.
Because of the non-convexity of the 0-norm, we minimize the $\ell_1$-norm as a convex
relaxation:
\begin{equation}
 \label{Lone}
 \min \|W^i \|_1\quad \mathrm{subject}\;\mathrm{to}\; W^i\,M = R^i,\;\; W^i \geq 0,
\end{equation}
which is in the form of linear programming.
The fact that the data may in general contain noise suggests
us to solve the following unconstrained optimization problem,
\begin{equation}
 \label{LoneU}
 \min_{W^i \geq 0} \mu\|W^i \|_1 + \frac{1}{2}\|R^i - W^i\,M \|^2_2\;,
\end{equation}
for which Bregman iterative method \cite{G_O_09, YO} with a proper projection onto non-negative convex subset
is used to obtain a solution.  In this paper, we shall use the linearized Bregman
iteration to solve (\ref{LoneU}) due to its efficiency.  For each row $W^i$ of $W$, we introduce $u = (W^i)^{\mathrm{T}}, f = (R^i)^{\mathrm{T}}, B = M^{\mathrm{T}}$, then (\ref{LoneU}) is equivalent to
\begin{equation}
\label{linBreg}
\min_{u \geq 0} \mu\|u \|_1 + \frac{1}{2}\|f - B\,u \|^2_2\;.
\end{equation}
The $l_2$ norm in (\ref{LoneU},\ref{linBreg}) is to model the unknown measurement error or noise as Gaussian.
When there is minimal measurement error, one must assign a tiny value to $\mu$ to heavily weigh the fidelity term $\displaystyle \|f - B\,u \|^2_2$ in order for $\displaystyle B\,u = f $ to be nearly satisfied.
The linearized Bregman method can be written iteratively by introducing an
auxiliary variable $v^j$:
\begin{equation}
\left\{
\begin{array}{ll}
 v^{j + 1} = v^j - B^\mathrm{T}(B\,u^j - f),\\
 \vspace {1.5mm}
 u^{j + 1} = \delta \cdot shrink_{+}  (v^{j + 1}, \mu \bigr),\,\,
       \end{array}
\right.
\end{equation}
where $u^0 = v^0 = O $, $\delta > 0 $ is the step size, and $shrink_{+}$ is for computing nonnegative solutions,
\begin{equation}
shrink_{+}  (v, \mu \bigr) =
\left\{
\begin{array}{ll}
v - \mu , \hskip 0.5 mm \quad \mathrm{if}\, v>\mu\,,\\
0, \qquad\quad \mathrm{if}\, v<\mu\,.
\end{array}
\right.
\end{equation}

\section{Numerical experiments}
In this section, we shall apply our method to the so-called SWOrRD Raman data provided by the Naval Research Laboratory, and present the computational results.
\subsection{SWOrRD Data}
\begin{figure}
\includegraphics[height=8cm,width=9cm]{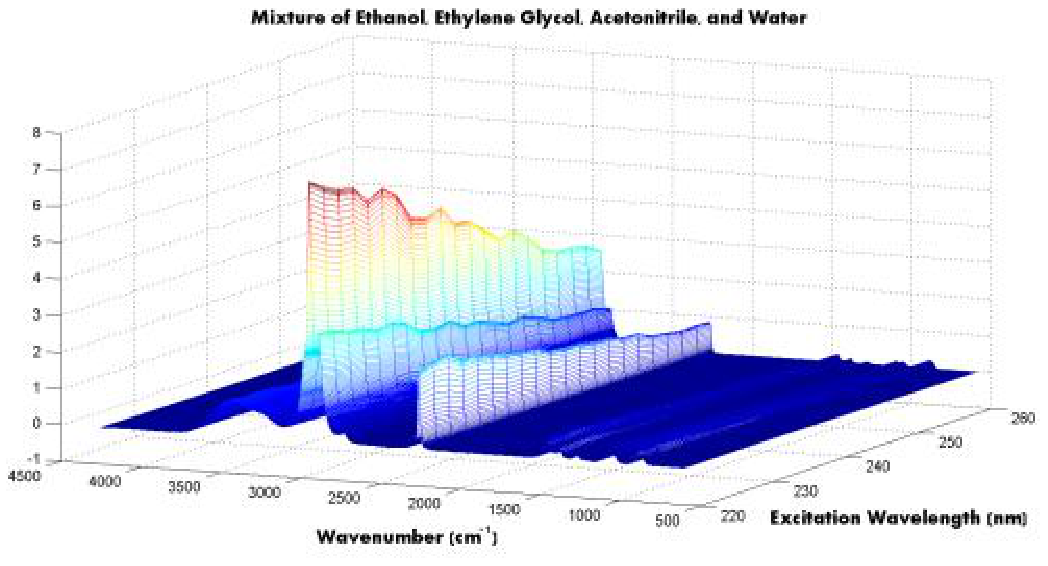}
\includegraphics[height=8cm,width=8cm]{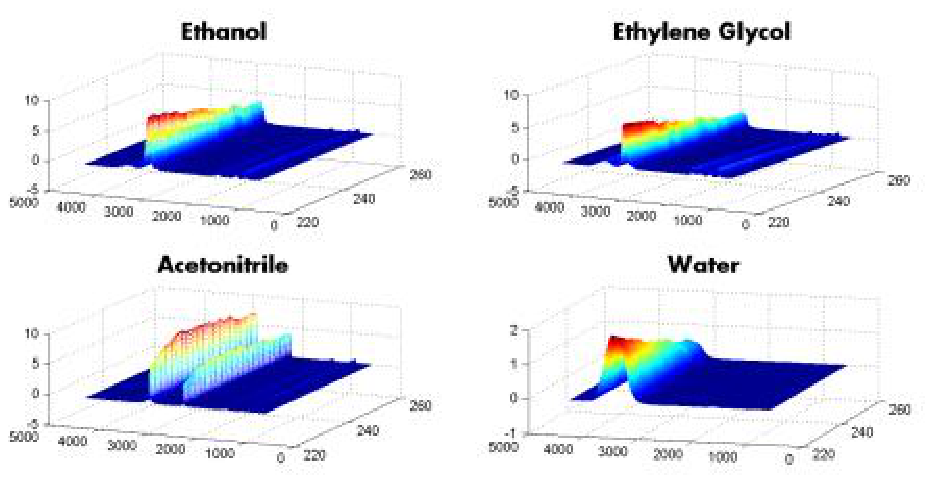}
\caption{SWorRD 2D signatures of four pure substances are shown in boxes on the right.  A signature of a mixture of these pure
substances is shown on the left. }
 \label{raman_mix4}
\end{figure}

The Swept Wavelength Optical Resonant Raman Detector (SWOrRD) at the Naval Research Laboratory is
a spectroscopy system which is able to produce two dimensional spectral data of biological agents
and chemical substances.  A target substance is illuminated with a specific laser wavelength and this
generates a resonance Raman spectrum.   When the laser wavelength is varied and the
process is repeated, the resonance Raman spectra (one at each laser wavelength),
forms a two-dimensional plot (signature) where
one axis is the input laser wavelength and the other axis is the
wavenumber of the Raman spectrum.
Fig. \ref{raman_mix4} is an example of 2D SWorRD Raman spectra
of a mixture of ethanol, ethylene glycol, acetonitrile, and water.
The data matrix $X$ corresponds to this 2D spectra.
Each column $X$ is a Raman spectrum of the mixture at a specific laser wavelength.
For SWorRD data, there are much more mixtures than the number of sources.
The spectral features may be slightly altered from one
wavelength to another.  To reduce this effect, we do not necessarily use all the data,
and so we sample the spectra at a subset of input laser wavelengths.
\subsection{Results}
\begin{figure}
\includegraphics[height=6cm,width=8cm]{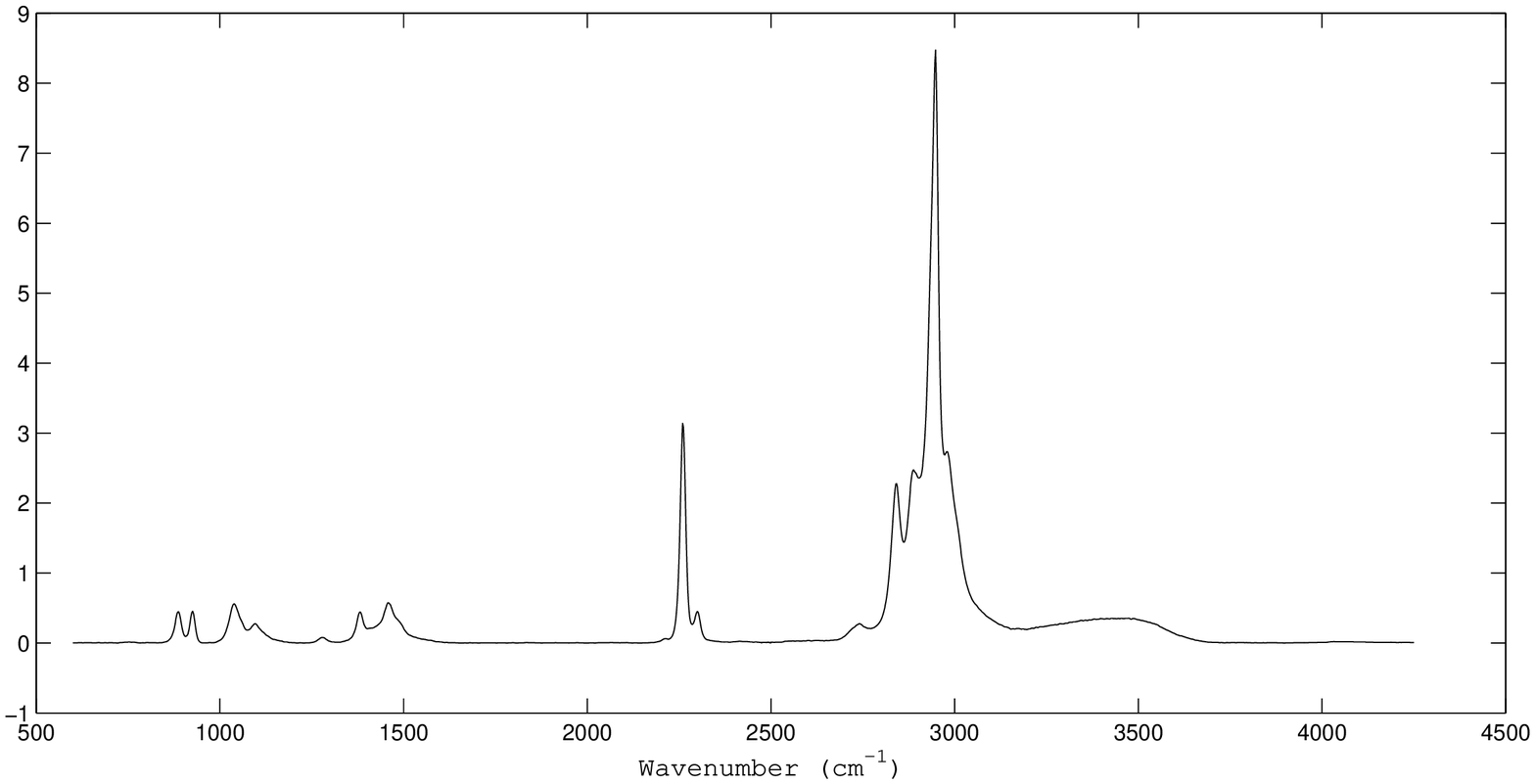}
\includegraphics[height=6cm,width=8cm]{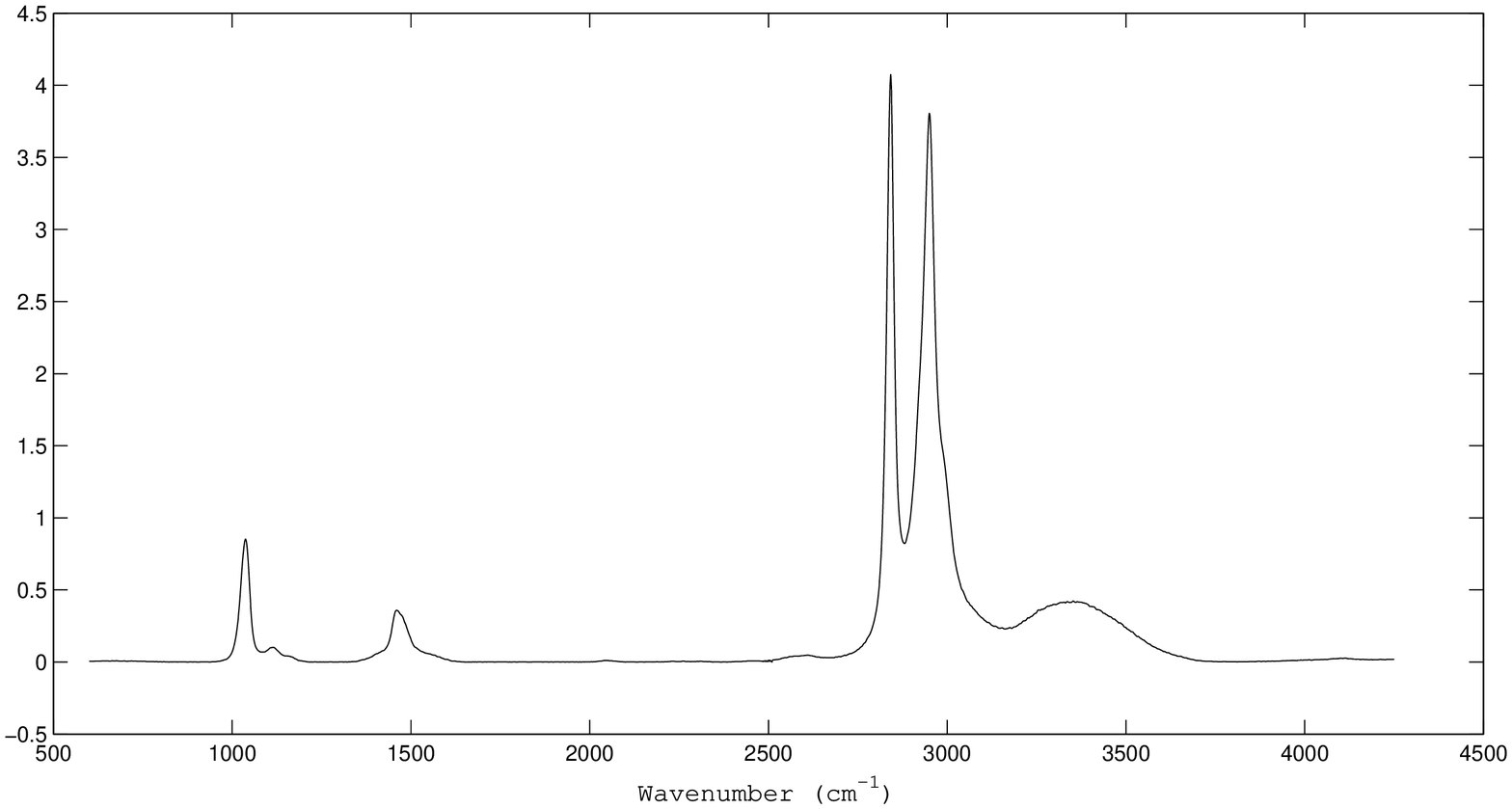}
\caption{A Raman spectrum of the mixture of methanol and two other liquid chemicals is shown on the left, while the
Raman spectral reference of methanol is plotted on the right.}
\label{result1}
\end{figure}

\begin{figure}
\includegraphics[height=8cm,width=16cm]{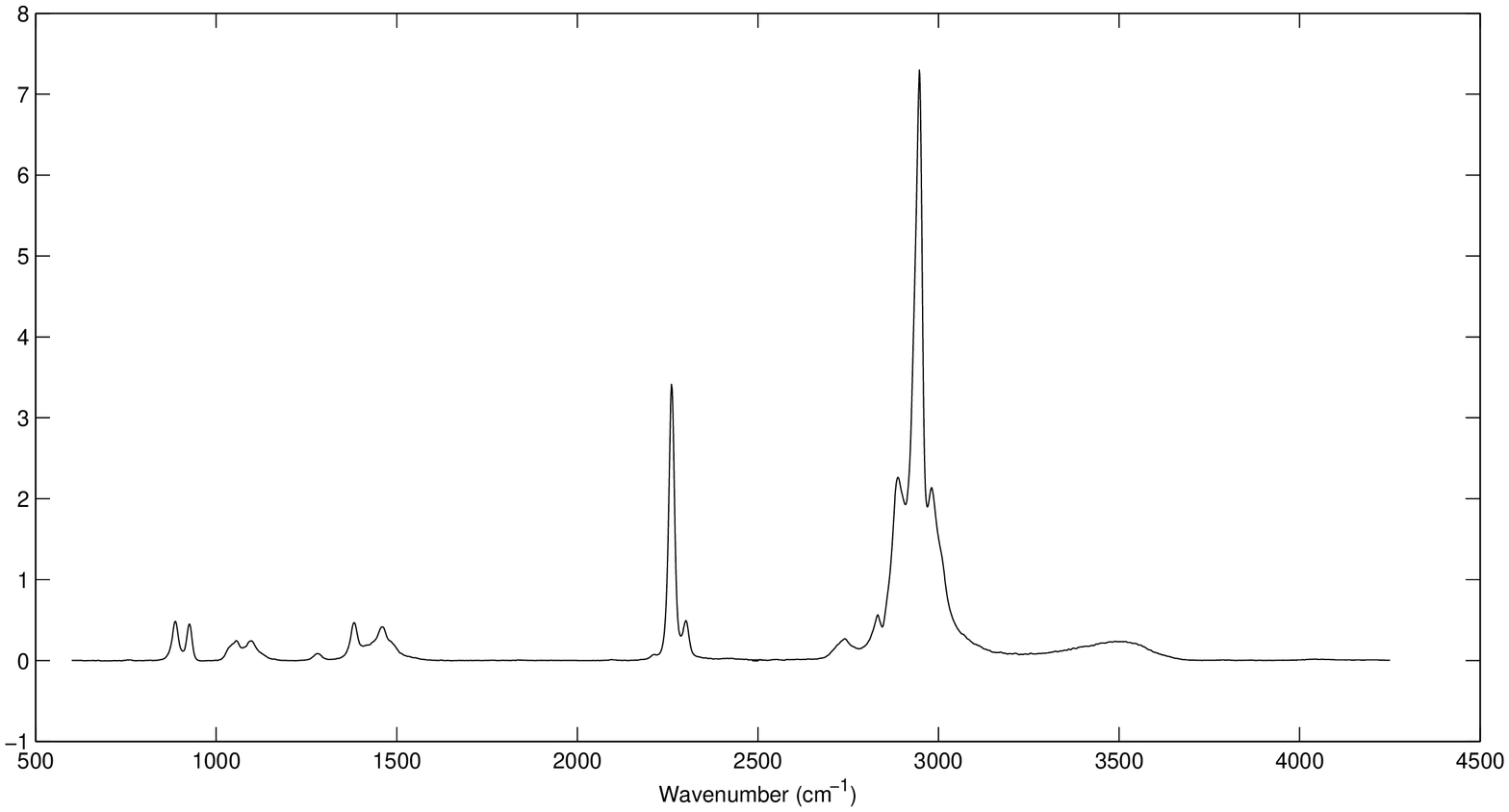}
\caption{One fitting residual from the constrained least squares data fitting.}
\label{result2}
\end{figure}

\begin{figure}
\includegraphics[height=6cm,width=8cm]{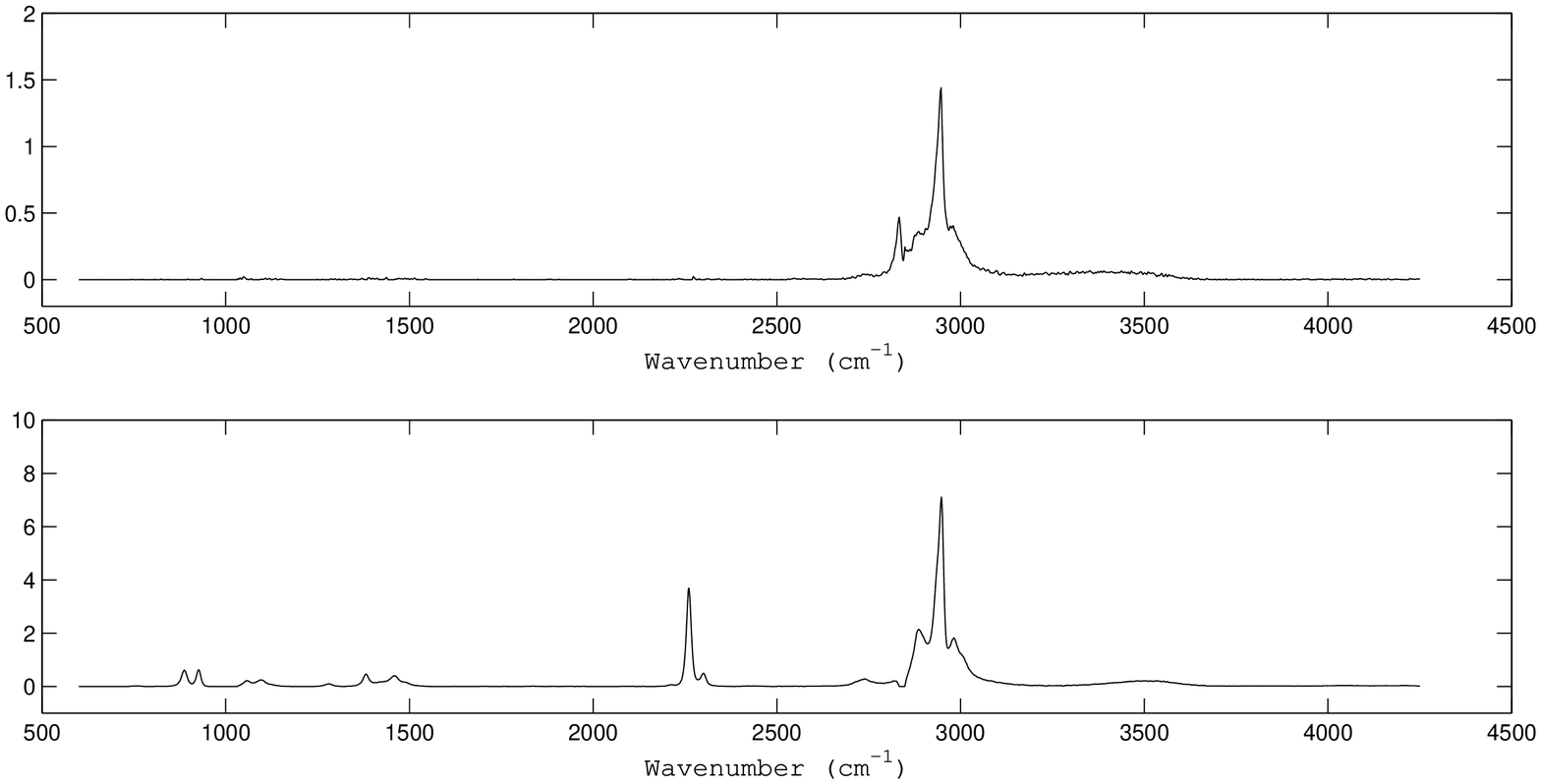}
\includegraphics[height=6cm,width=8cm]{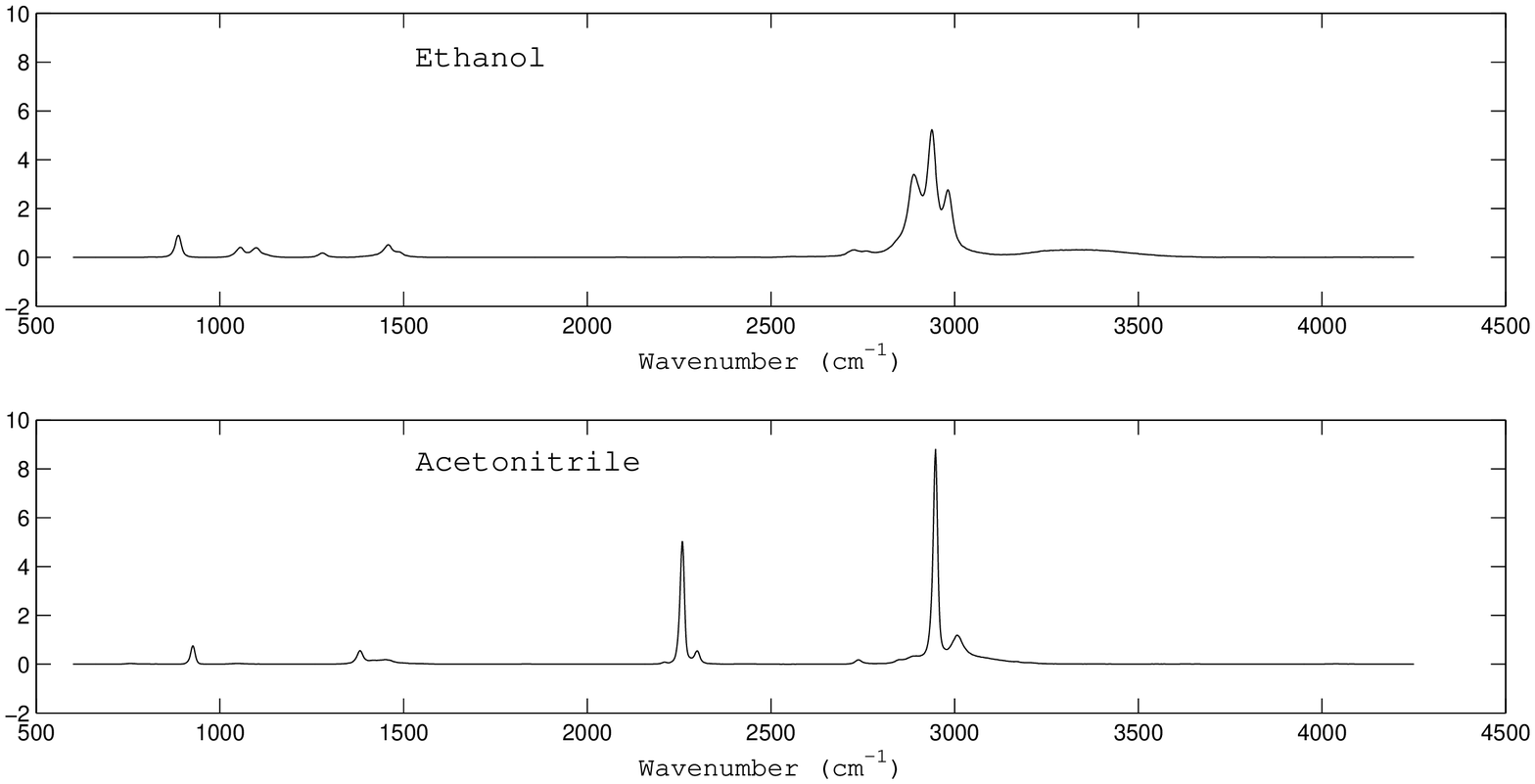}
\caption{Left column is the two identified spectral structures, right column shows the Raman spectral references for
ethanol and acetonitrile.  The value of $\mu$ used in step 2 is 0.09.}
\label{result3}
\end{figure}

\begin{figure}
\includegraphics[height=6cm,width=8cm]{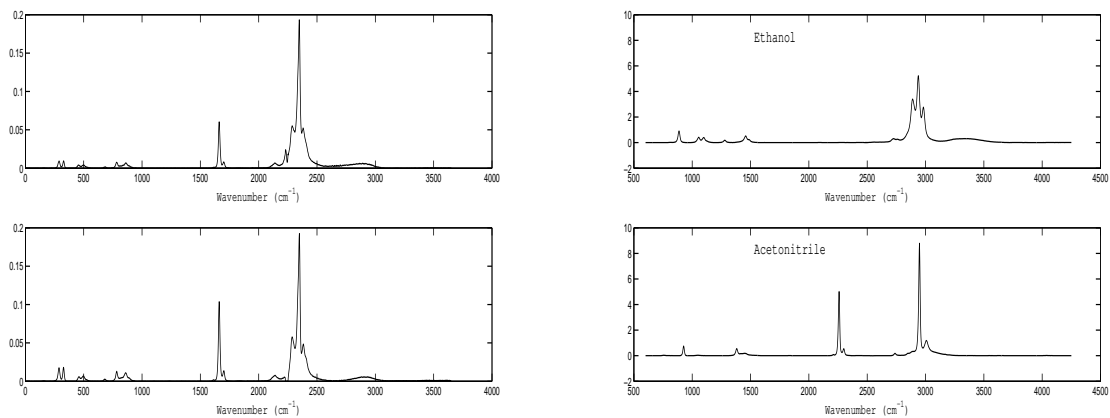}
\includegraphics[height=6cm,width=8cm]{figures/ref_Eth_Ace.eps}
\caption{Left column is the results computed by NMF, right column shows the Raman spectral references for
ethanol and acetonitrile.}
\label{result4}
\end{figure}

\begin{figure}
\includegraphics[height=6cm,width=8cm]{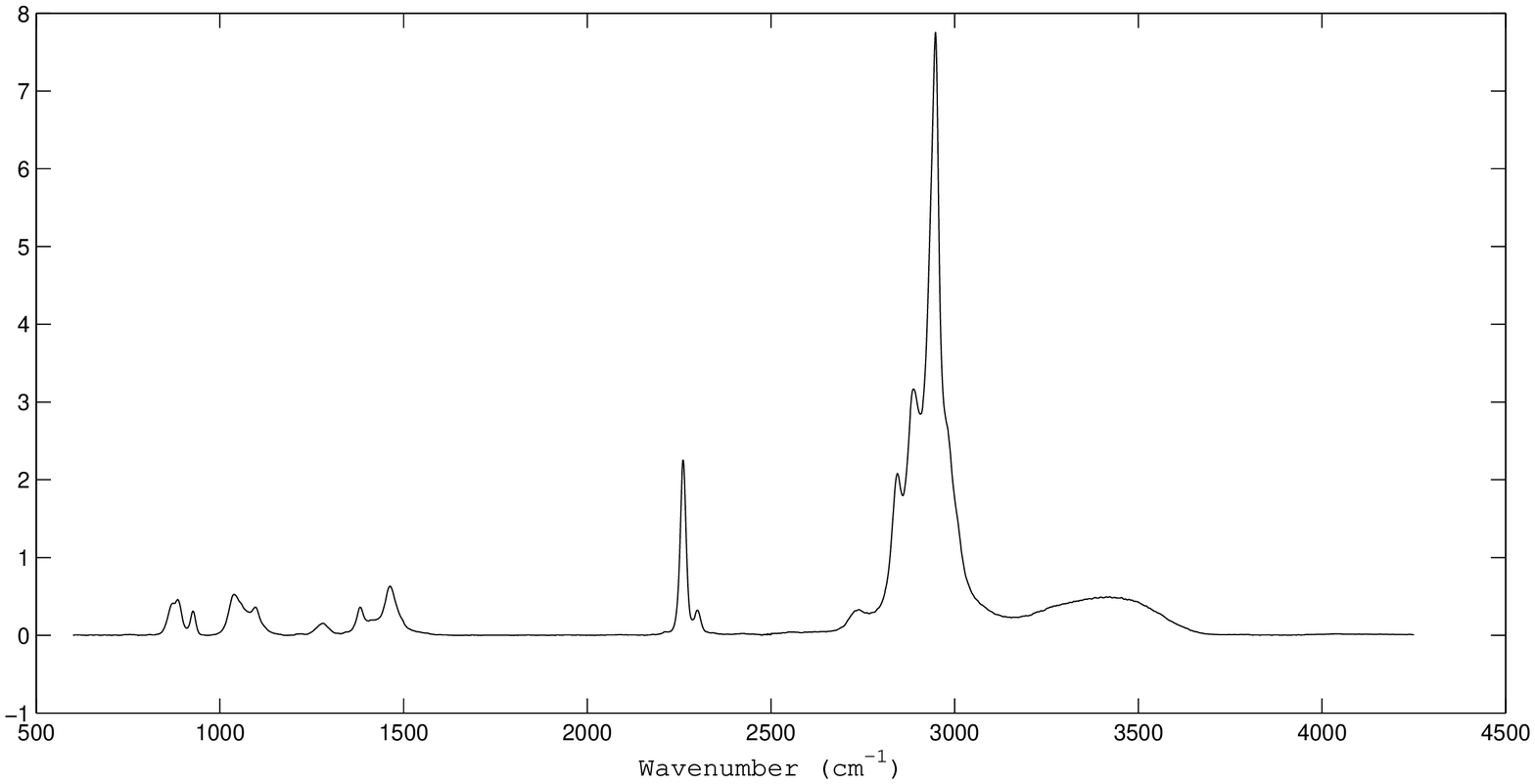}
\includegraphics[height=6cm,width=8cm]{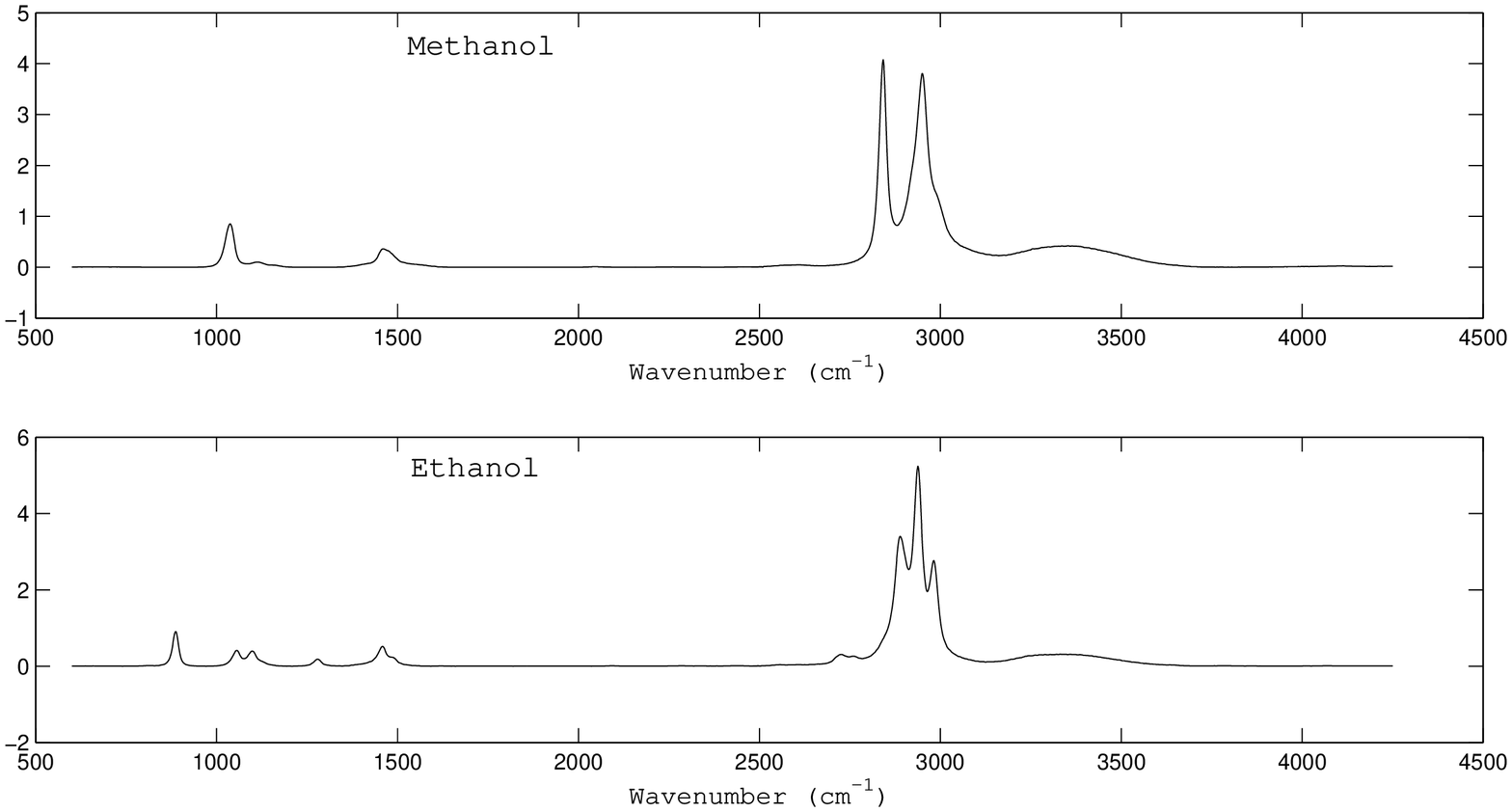}
\caption{A Raman spectrum of the mixture of methanol, ethanol and two other chemicals is shown on the left, while the
Raman spectral references of methanol and ethanol are plotted on the right.  $\mu = 0.09$ is used for the source recovery in step 2. }
\label{eg2_1}
\end{figure}

\begin{figure}
\includegraphics[height=8cm,width=16cm]{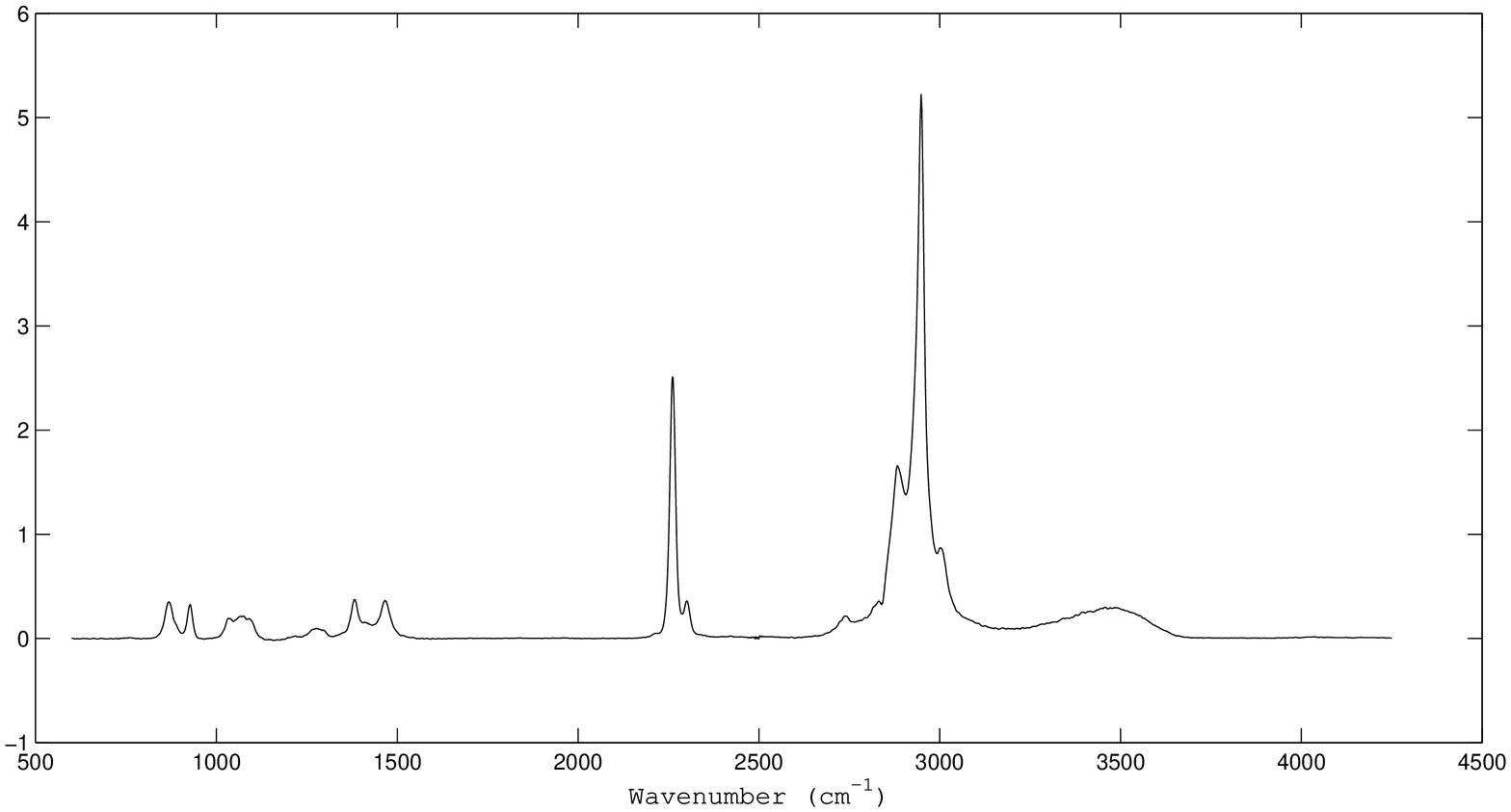}
\caption{One fitting residual from the constrained least squares data fitting.}
\label{eg2_2}
\end{figure}

\begin{figure}
\includegraphics[height=6cm,width=8cm]{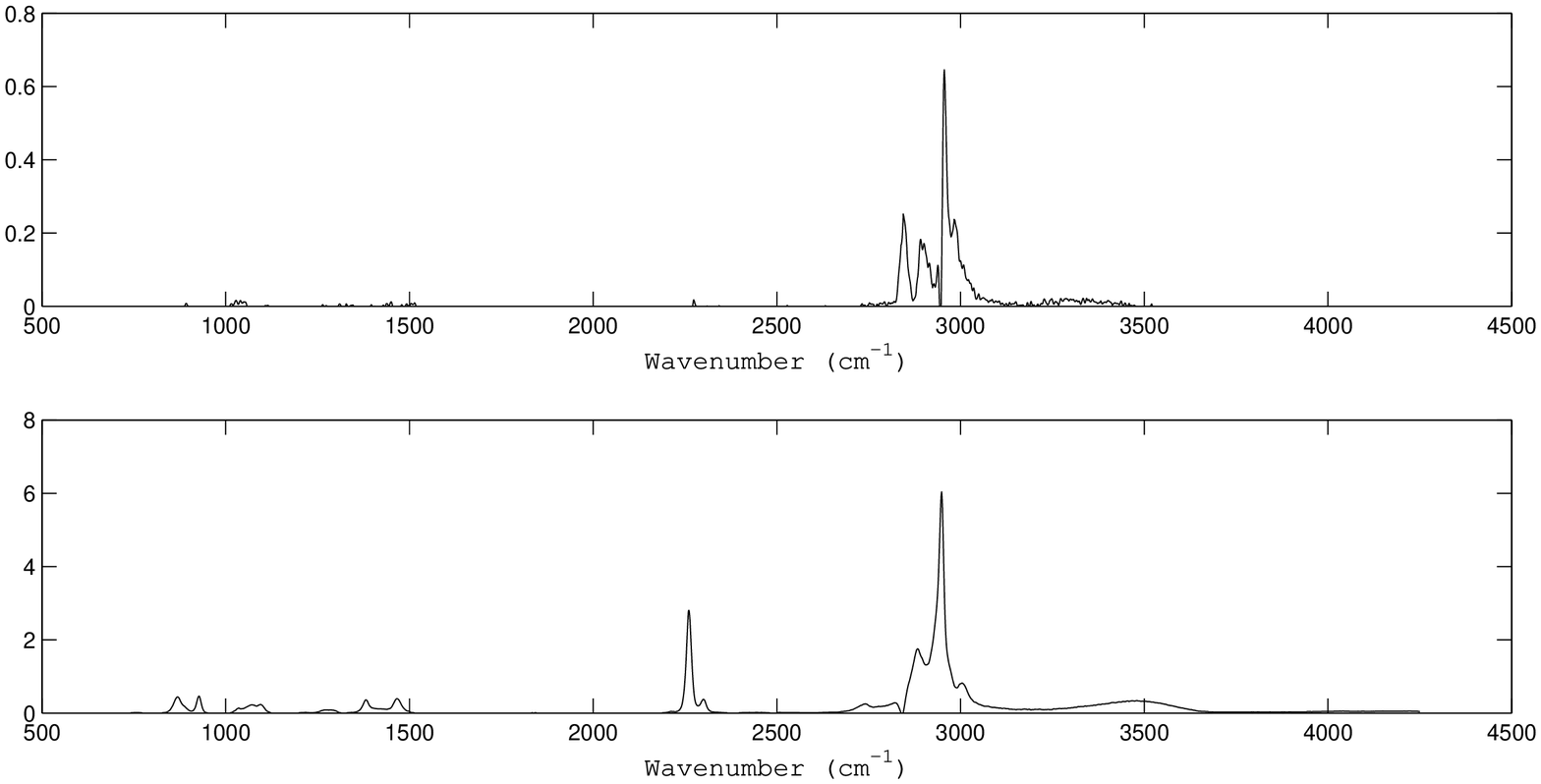}
\includegraphics[height=6cm,width=8cm]{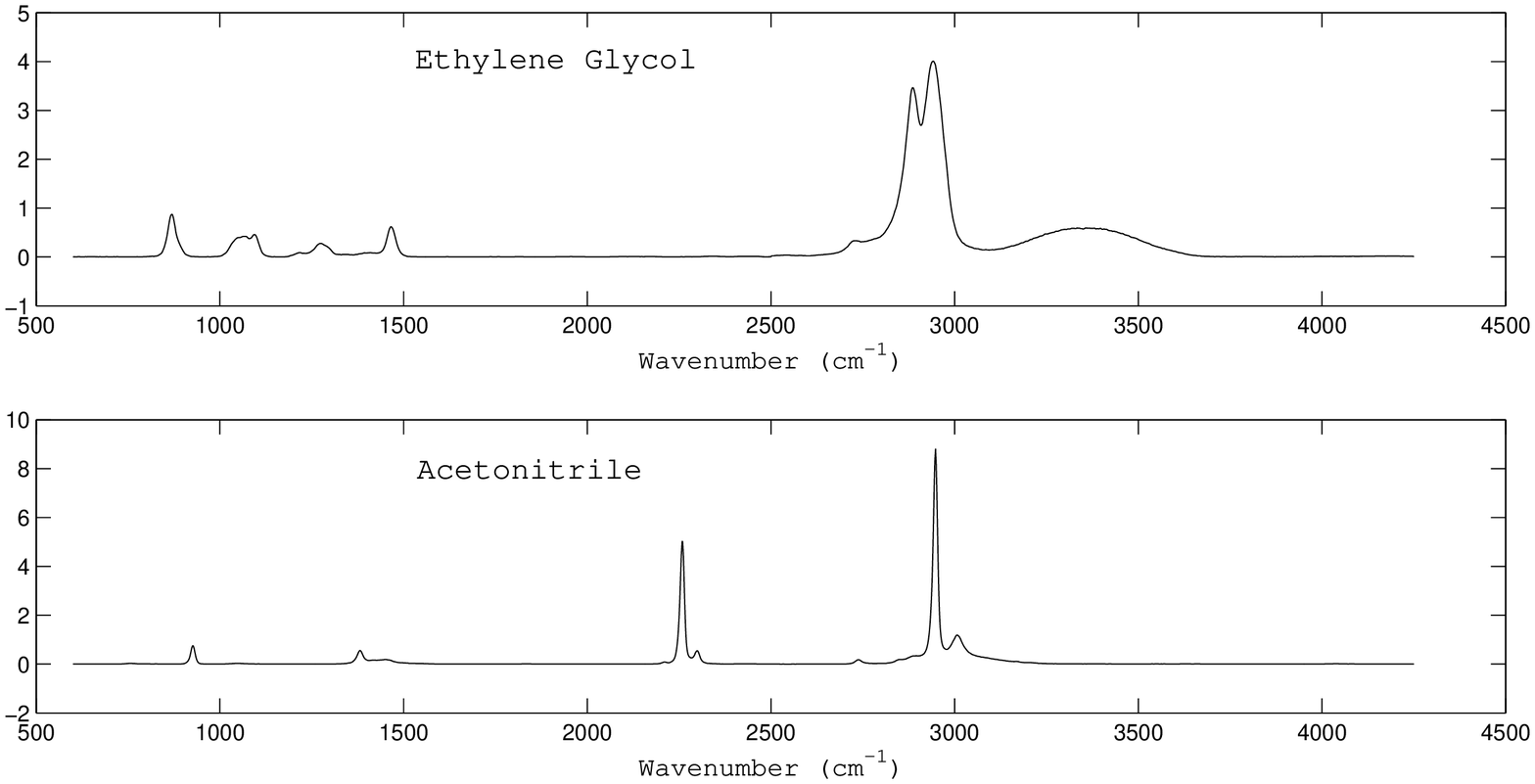}
\caption{Left column is the two identified spectral structures, right column shows the Raman spectral references for
acetonitrile and ethylene glycol.}
\label{eg2_3}
\end{figure}

\begin{figure}
\includegraphics[height=8cm,width=16cm]{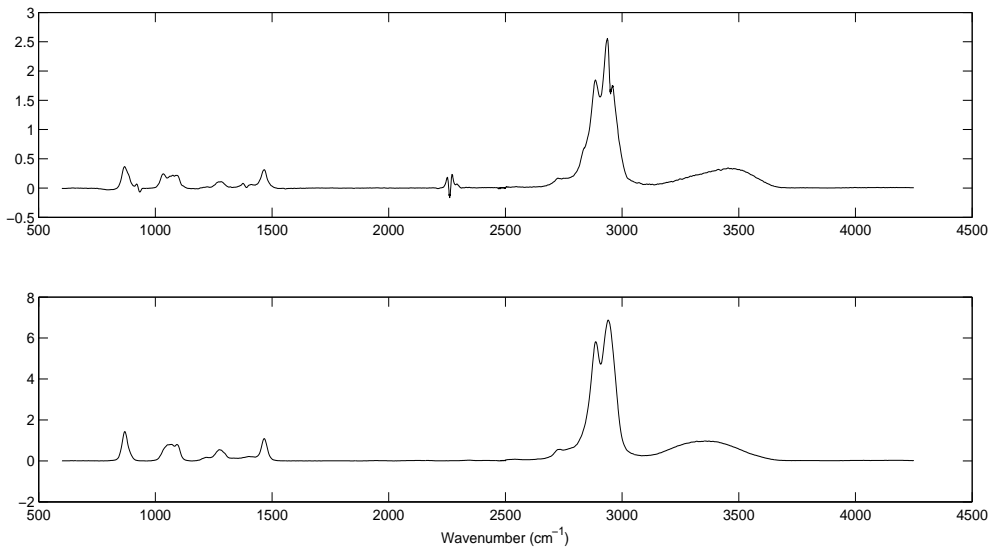}
\caption{Top panel is the recovered structure after subtracting acetonitrile.  The spectral line shape is easily recognizable as ethylene glycol, whose spectral reference is shown in the bottom plot. }
\label{eg2_4}
\end{figure}

We report here the computational results from our proposed method.
We tested our method on two sets of SWorRD data.
The samples consist of several liquid substances, some of them are commonly
used for synthesizing explosives.
The first data set includes 21 mixed Raman spectra at different excitation wavelengths.
This sample is known to contain methanol and its concentration is below $1/3$.
There are also two {\em a-priori} unknown liquid chemicals.
The method first fits the known reference spectrum of methanol to the Raman data.
Then it identifies two unknown chemicals from the fitting residuals.
The results are shown in a series of plots, Fig. \ref{result1}-Fig. \ref{result4}.
From all the data available, we found that the 5 spectra at consecutive laser
wavelengths of $(248,250,252,254,256)$ nm produce the best results.
The theoretical underpinning of optimal selection of the data is under further study.
 Fig. \ref{result1} shows the line shape of the spectrum of the mixture, and
the spectral reference of methanol.  The residual after fitting methanol
to the data is plotted in Fig.  \ref{result2}, where some structure can be seen.
Then further identification of the two hidden chemicals was made by the convex BSS method,
and the results are presented in Fig. \ref{result3}.
Compared to the ground truth, the results are satisfactory in that
the recovered spectral structures are recognizable as ethanol and acetonitrile.
Conventional method such as NMF tends to lose reliability in dealing with the data,
the NMF results are also presented as a comparison.  In the second example, we try to identify two hidden chemicals from a mixture of four chemicals. The known chemicals are methanol and ethanol. The total concentration of the two
is below $1/2$, which is known from the sensing hardware.
The computational results are presented in Fig. \ref{eg2_1}-Fig. \ref{eg2_4}.
Although no apparent Raman spectral structure can be identified from the first plot in Fig. \ref{eg2_3} ,
the second structure is easily recognizable as acetonitrile upon comparison with reference spectrum.
Next we subtract the confirmed acetonitrile from the residual by solving the constrained least squares in step 1,
then feed the new residual back to step 2 for extracting more hidden chemicals.
The extracted structure in Fig. \ref{eg2_4} can be easily identified as ethylene
glycol comparing to the ground truth.  Computations are performed on a Dell laptop with 6G RAM and 1.6 GHz i7 CPU.
The cpu time for example 1 is 7.613 seconds, and it is 10.327 seconds for example 2.  The method proves to be efficient and can be of use in rapid detection of chemical substances.




\section{Conclusions and Future Work}
A semi-blind sparse and convex source
identification method is developed
to extract meaningful spectral structures from mixture data of Raman spectroscopy.
The method is designed to identify potentially hidden chemicals
after fitting the known spectral references to the data.  Our method can be useful
for extracting unknown source signals from the residuals
after known reference spectra have been first deployed to fit the data.
The major strength of the technique is its ability
to be used either with known reference spectra for quantification or
without reference spectra for identification of unknown/hidden chemical substances.
Numerical results on SWOrRD data showed the promising potential of
our method on explosives detection.

The model considered in the paper
is a linear and stationary model which assumes no shift/squeeze in the spectral lines.
A future line of work will study how to build this nonlinear effect into
the identification model.  Given that the shift/squeeze amount is small,
one may explore the idea of image registration.
We also plan to study more reliable and efficient methods for
the residuals decomposition, as their success highly depends
on a viable working assumption on the pure signals.

The semi-blind source identification problem we addressed here
also has analogues in detecting atmospheric trace gases with the so called differential
optical absorption spectroscopy.
Analysis of the fitting residuals for mysterious species arises there as well, and
is non-convex in general. Recently the authors achieved some success
in this direction based on a similar semi-blind framework,
although the details of fitting and blind source identification process are quite
different \cite{Sun_Xin_Doas}.

\end{document}